\documentclass[12pt]{article}
\usepackage{amssymb}
\usepackage{hyperref}
\usepackage{amsmath,amscd}
\hypersetup{pdfpagemode=FullScreen,colorlinks=true}

\newtheorem{pro}{Proposition}
\newtheorem{dfi}{Definition}
\newtheorem{ntn}{Notation}
\newtheorem{thm}{Theorem}
\newtheorem{cor}{Corollary}
\newtheorem{rem}{Remark}
\newtheorem{lem}{Lemma}
\newtheorem{exa}{Example}

\def\R{\mathbb{R}}
\def\H{\mathbb{H}}
\newcommand{\he}[1]{{\mathbb H}^{#1}}

\newenvironment{pf}{\begin{trivlist}\item[]{\bf Proof\ }}
{\mbox{}\hfill\rule{.08in}{.08in}\end{trivlist}}

\title{Cup-products in $L^{q,p}$-cohomology: discretization and quasi-isometry invariance}
\author{Pierre Pansu\footnote{P. P.~is supported by MAnET Marie Curie Initial Training Network and Agence Nationale de la Recherche grants ANR-2010-BLAN-116-01 GGAA and ANR-15-CE40-0018 SRGI. He gratefully acknowledges the hospitality of Isaac Newton Institute, where this work was completed, of EPSRC and of Simons Foundation.
}}

\begin{document}

\maketitle

\abstract{We relate $L^{q,p}$-cohomology of bounded geometry Riemannian manifolds to a purely metric space notion of $\ell^{q,p}$-cohomology, \emph{packing cohomology}. This implies quasi-isometry invariance of $L^{q,p}$-cohomology together with its multiplicative structure. The result partially extends to the Rumin $L^{q,p}$-cohomology of bounded geometry contact manifolds.}

\tableofcontents

\section{Introduction}

$L^{q,p}$ cohomology is a quantitative variant of de Rham cohomology for Riemannian manifolds: differential forms are required to belong to $L^p$ spaces, i.e. to decay at infinity. It has proven its usefulness in various parts of geometry and topology, \cite{GKahler}, \cite{Luck}, \cite{BP}, \cite{BK}.

Because of its topological origin, it is expected that $L^{q,p}$ cohomology be computable by many different means, and be quasi-isometry invariant. This has been established over the years in many cases, \cite{GAs}, \cite{Fan}, \cite{Elek}, \cite{BP}, \cite{Ducret}, \cite{Genton}. In this paper, one completes the picture, 
\begin{itemize}
  \item by covering all remaining cases (limiting cases for exponents $p$ and $q$),
  \item by proving invariance of cup-products.
\end{itemize}

The new input is two-fold. 
\begin{enumerate}
  \item We exploit progress made in the 2000's on the $L^1$ analytic properties of the exterior differential, \cite{BB2}, \cite{VS}, \cite{LS}.
  \item We use a definition of $\ell^{q,p}$ cohomology for metric spaces, \emph{packing cohomology}, which is well-suited to handle products.
\end{enumerate}

Next we proceed to precise statements. Packing cohomology will be defined in subsection \ref{pack}.

\begin{dfi}
Let $X$ be a Riemannian $n$-manifold. $\Omega^{q,p,k}(X)$ denotes the space of $L^{q}$ differential forms whose distributional exterior derivative is an $L^{p}$ differential form. Define $L^{q,p}$ cohomology by
\begin{eqnarray*}
H^{q,p,k}(X)=\mathrm{ker}(d)\cap \Omega^{p,p,k}(X)/d\Omega^{q,p,k-1}(X).
\end{eqnarray*}
\emph{Exact $L^{q,p}$-cohomology} $EH^{q,p,k}(X)$ is the kernel of the forgetful map $H^{q,p,k}(X)\to H^k(X,\R)$.
\end{dfi}

\begin{dfi}
Say a metric space $X$ has uniformly vanishing cohomology up to degree $k$ if, for every $R$, there exists $\tilde{R}$ such that for every $j=0,\ldots,k$ and every $x\in X$, the map $H^j(B(x,\tilde{R}),\R)\to H^j(B(x,R),\R)$ induced by inclusion $B(x,R)\subset B(x,\tilde{R})$ vanishes.
\end{dfi}

\begin{thm}
\label{qib}
Assume that $1\leq p\leq q\leq \infty$ and $\frac{1}{p}-\frac{1}{q}\leq \frac{1}{n}$. Consider the class of Riemannian manifolds with the following properties.
\begin{enumerate}
\item Dimension equals $n$.
  \item Bounded geometry: there exist $M>0$ and $r_0>0$ and for every point $x$ an $M$-bi-Lipschitz homeomorphism of the unit ball of $\R^n$ onto an open set containing $B(x,r_0)$.
  \item Uniform vanishing of cohomology up to degree $k-1$.
\end{enumerate}
If $p=1$, $q=n/n-1$ and $k=n$, one should replace $L^{n/n-1,1}$-cohomology with $L^{n/n-1,\mathcal{H}^1}$-cohomology, to be defined in subsection \ref{lpi}. If $p=n$, $q=\infty$ and $k=1$, one should replace $L^{\infty,n}$-cohomology with $L^{BMO,n}$-cohomology, to be defined in subsection \ref{lpi} as well.

For $X$ in this class, and up to degree $k-1$, $L^{q,p}$-cohomology and packing $\ell^{q,p}$-cohomology of $X$ at all sizes are isomorphic as vectorspaces. Furthermore, in degree $k$, the exact $L^{q,p}$-cohomology and exact packing $\ell^{q,p}$-cohomology of $X$ at all sizes are isomorphic. Finally, these spaces, together with their multiplicative structure, are quasi-isometry invariant.

If $(p,q)\not=(1,\frac{n}{n-1})$, or $(n,\infty)$, the isomorphisms are topological, they arise from homotopy equivalences of complexes of Banach spaces.
\end{thm}

Along the way, we shall establish an analogue of Theorem \ref{qib} (except its multiplicative content) for contact manifolds equipped with bounded geometry Carnot-Carath\'eodory metrics and the Rumin complex. This relies on recent $L^1$ analytic results for invariant operators on Heisenberg groups, \cite{CVS}, \cite{SmP}. It would be nice to extend this result to larger classes of equiregular Carnot manifolds. The machinery developped here would yield it provided the needed analytical properties of Rumin's complex were known. Unfortunately, Rumin's complex does not form a differential algebra, so it cannot capture the multiplicative structure of cohomology.

\subsection{Plan of the paper}

Section \ref{input} collects the needed Euclidean Poincar\'e inequalities. Section \ref{leray act} recalls Leray's proof of de Rham's theorem relating de Rham to \v Cech cohomology. Section \ref{leray cust} presents a new variant of Leray's method, which is far less demanding in terms of properties of coverings and Poincar\'e inequalities. The loss on domains in Poincar\'e inequalities that it allows is crucial in two ways,
\begin{enumerate}
  \item It feeds on existing, perhaps suboptimal in terms of domains, analytical inequalities.
  \item It allows to jump from one scale to a much larger one, under a mere global topological assumption.
\end{enumerate}
This is illustrated in section \ref{leray simpl}, where the $\ell^p$ cohomology of a simplicial complex with uniformly vanishing cohomology is shown to coincide with that of its Rips complex at arbitrary scales. In section \ref{qi inv}, this result is reformulated in terms of Alexander-Spanier cochains and packing cohomology, a theory which is quasi-isometry invariant by nature. Note that the main output of sections \ref{leray simpl} and \ref{qi inv} (functoriality of $\ell^{q,p}$-cohomology of simplicial complexes under coarse embeddings) is valid with no other restriction on $(q,p)$ than $1\leq p\leq q\leq +\infty$).
Section \ref{contact} details the analogous result for contact sub-Riemannian manifolds. Some extra analytical difficulties arise since the adapted exterior differential, due to M. Rumin, is a second order operator in middle dimension.

\section{Analytical input}
\label{input}

\subsection{Poincar\'e inequalities}
\label{poincar}

We shall use the following results, which can be found in \cite{SmP}.

\begin{thm}[Baldi-Franchi-Pansu]
\label{poincare}
Assume that $1\leq p\leq q\leq \infty$ and $\displaystyle \frac{1}{p}-\frac{1}{q}\leq \frac{1}{n}$. Let $\lambda>1$. Let $B=B(0,1)$ and $B'=B(0,\lambda)$ be concentric balls of $\R^n$. 

Assume first that $(p,q,k)\notin\{(1,n/n-1,n),(n,\infty,1)\}$. There exists a constant $C=C(\lambda)$ such that for every closed differential $k$-form $\omega$ on $B'$, there exists a differential $k-1$-form $\phi$ on $B$ such that $d\phi=\omega_{|B}$ and
\begin{eqnarray*}
\|\phi\|_{L^{q}(B)}\leq C\,\|\omega\|_{L^{p}(B')}.\quad \quad \quad (Poincare_{q,p}(k))
\end{eqnarray*}

If $p=1$, $q=n/n-1$ and $k=n$, inequality $(Poincare_{q,p}(k))$ is replaced with
\begin{eqnarray*}
\|\phi\|_{L^{n/n-1}(B)}\leq C\,\|\omega\|_{\mathcal{H}^{1}(B')}.\quad \quad \quad (Poincare_{n/n-1,\mathcal{H}^1}(n))
\end{eqnarray*}

If $p=n$, $q=\infty$ and $k=1$, inequality $(Poincare_{q,p}(k))$ is replaced with
\begin{eqnarray*}
\|\phi\|_{BMO(B)}\leq C\,\|\omega\|_{L^{n}(B')}.\quad \quad \quad (Poincare_{BMO,n}(1))
\end{eqnarray*}
\end{thm}
Similar inequalities hold for $\lambda$ large enough on Heisenberg balls, with exterior differential $d$ replaced with Rumin's differential $d_c$, see \ref{poincaredc}.

\begin{rem}
Note that inequalities $(Poincare_{n/n-1,1}(n))$ and \break $(Poincare_{\infty,n}(1))$ fail.
\end{rem}

\subsection{$L^{BMO}$ and $L^{\mathcal{H}^1}$ norms}
\label{lpi}

To cover the exceptional configurations $p=1$, $q=n/n-1$ and $k=n$, on one hand, and $p=n$, $q=\infty$ and $k=1$, on the other hand, one needs switch from Lebesgue spaces to mixed Lebesgue-Hardy spaces.

\begin{dfi}
Let $X$ be a bounded geometry Riemannian manifold. For a differential forms $\omega$ on $X$, define
\begin{eqnarray*}
\|\omega\|_{L^{BMO}}=\sup_{x\in X}\|\omega\|_{BMO(B(x,1))}, \quad \|\omega\|_{L^{\mathcal{H}^1}}=\int_{X}\|\omega\|_{\mathcal{H}^1(B(x,1))}\,dx. 
\end{eqnarray*}
\end{dfi}
These are the norms used in the definition of $L^{BMO,n}$ and $L^{n/n-1,\mathcal{H}^1}$-cohomology, required only in degrees $1$ and $n$ respectively. One does not need modify the definition of packing $\ell^{\infty,n}$ and $\ell^{n/n-1,1}$-cohomology.

\subsection{$L^\pi$-cohomology}

For the proofs, it will be necessary to deal with a whole complex at the same time. 

\begin{ntn}
Let $\pi=(p_0,\cdots,p_n)$, where $1\leq p_k\leq\infty$ for $k=0,\ldots,n$. $\Omega^{\pi,k}(X)$ denotes the space of $L^{p_k}$ differential forms whose distributional exterior derivative is an $L^{p_{k+1}}$ differential form. The norm there is 
$$|\omega|_{p_k}+|d\omega|_{p_{k+1}}.$$ 
The exterior differential $d$ is a bounded operator on 
$$
\Omega^{\pi,\cdot}(X):=\bigoplus_{k=0}^n \Omega^{\pi,k}(X).
$$
It constitutes a complex whose cohomology 
\begin{eqnarray*}
H^{\pi,\cdot}(X)=\mathrm{ker}(d)\cap \Omega^{\pi,\cdot}(X)/d\Omega^{\pi,\cdot}(X)
\end{eqnarray*}
is called the \emph{$L^{\pi}$-cohomology} of $X$. \emph{Reduced $L^{\pi}$-cohomology} $ \bar{H}^{\pi,\cdot}(X)$ is obtained by modding out by the closure of the image of $d$.
\end{ntn}

Note that, for $k=0,\ldots,n$, $H^{p_{k-1},p_{k},k}(X)=H^{\pi,k}(X)$ for any sequence $\pi$ containing $(p_{k-1},p_{k})$ as a subsequence.

\section{Leray's acyclic coverings theorem}
\label{leray act}

Let $X$ be a Riemannian manifold. Let $\mathcal{U}=(U_i)_{i\in I}$ be an open coverings of $X$. Assuming that Poincar\'e's inequality holds as in Proposition \ref{poincare} for all pairs $(U_i,U_{i})$ and all intersections $U_{i_1\ldots i_k}:=U_{i_1}\cap\cdots\cap U_{i_k}$ with uniform constants, we shall show that $L^{\pi}$-cohomology of $X$ is isomorphic to the $L^{\pi}$-cohomology of the \emph{nerve} $T$, of $\mathcal{U}$, i.e. the simplicial complex which has a vertex $i$ for each open set $U_i$ and a face $i_1\ldots i_k$ each time the intersection $U_{i_1\ldots i_k}\not=\emptyset$. We shall furthermore assume that the nerve is locally bounded (every $U_i$ intersects a bounded number of other $U_j$'s), and we shall need a partition of unity $(\chi_i)$ subordinate to $\mathcal{U}$ such that the gradients $\nabla\chi_i$ are uniformly bounded.

Recall that simplicial cochains are skew-symmetric functions on oriented simplices.

\subsection{Closed 1-forms and 1-cocycles}

Let us first explain the argument for $L^{\pi}H^1(X)$. 

Given a closed 1-form $\omega$ on $X$, let us view the collection $\tilde{\omega}$ of its restrictions $\omega_i=\omega_{|U_i}$ as a 0-cochain with values in 1-forms, $\tilde{\omega}\in C^{0,1}$. It is a 0-cocycle,
\begin{eqnarray*}
\delta\tilde{\omega}_{ij}:={\omega_i}_{|U_{ij}}-{\omega_i}_{|U_{ij}}=0.
\end{eqnarray*}
By assumption, $d\tilde{\omega}=0$. Poincar\'e inequalities provide us, for each $U_i$, with a primitive $\phi_i$ of $\omega_i$, $d\phi_i=\omega_i$. This forms a 0-cochain $\tilde{\phi}=(\phi_i)_{i\in I}$ with values in 0-forms, $\tilde{\phi}\in C^{0,0}$. These 0-forms need not match on intersections, i.e. $\delta\tilde{\phi}_{ij}:={\phi_i}_{|U_{ij}}-{\phi_j}_{|U_{ij}}$ need not vanish. Note that $\tilde{\kappa}:=\delta\tilde{\phi}$ is a 1-cochain with values in 0-forms, $\tilde{\kappa}\in C^{1,0}$. Furthermore, $d\tilde{\kappa}=d\delta\tilde{\phi}=\delta d\tilde{\phi}=\delta\tilde{\omega}=0$. This means that each function $\tilde{\kappa}_{ij}$ is constant, one can view $\kappa$ as a real valued 1-cochain of the nerve. It is a cocycle, since $\delta\kappa=\delta\delta\tilde{\phi}=0$. 

Assume that $\omega\in L^{p_1}(X)$. Poincar\'e inequalities state that primitives $\phi_i$ have $L^{p_0}$-norms controlled by local $L^{p_1}$ norms of $\omega$, so $\tilde{\phi}\in\ell^{p_1}(L^{p_0})$. The coboundary is bounded, so $\tilde{\kappa}=\delta\tilde{\phi}\in\ell^{p_1}(L^{p_0})$. Since each $\kappa_i$ is constant, $\kappa\in\ell^{p_1}$. 

If different choices of primitives $\phi'_i$ are made, $\phi'_i=\phi_i+u_i$, then $u_i$'s form a 0-form valued cochain $\tilde{u}$ such that $\tilde{\phi}'=\tilde{\phi}+\tilde{u}$. Therefore $\tilde{\kappa}'=\tilde{\kappa}+\delta\tilde{u}$. Since $d\tilde{u}=0$, one can view $\tilde{u}$ as a real valued 1-cochain of the nerve, and 
$\kappa'=\kappa+\delta u$. Again, each (constant) $u_i$ is bounded by $\|\omega\|_{L^{p_1}(U_i)}$, so $u\in\ell^{p_1}$. Since $p_0\geq p_1$, $\ell^{p_1}\subset\ell^{p_0}$, thus $u\in\ell^{p_0}$. This shows that the cohomology class $[\kappa]\in L^{\pi}H^1(T)$ does not depend on choices.

If $\omega=d\phi$ is exact, with $\phi\in L^{p_0}(X)$, one can choose $\phi_i=\phi_{|U_i}$, $\delta\tilde{\phi}=0$, thus $\kappa=0$. Therefore we get a bounded linear map $L^{\pi}H^1(X)\to L^{\pi}H^1(T)$.

Conversely, given a 1-cocycle $\kappa$ of the nerve, i.e. a collection of real numbers $\kappa_{ij}$ such that $\delta\kappa=0$, we first view it as a 1-cocycle $\tilde{\kappa}$ with values in (constant) 0-forms. We set
\begin{eqnarray*}
\phi_i :=\sum_{\ell}\kappa_{i\ell}\chi_\ell.
\end{eqnarray*}
This defines a 0-form valued 0-cocycle $\tilde{\phi}\in C^{0,0}$. The map $\epsilon:C^{1,0}\to C^{0,0}$ that we just defined, is an inverse for $\delta$ on cocycles. Indeed, 
\begin{eqnarray*}
(\delta\tilde{\phi})_{ij}
&=&\sum_{\ell}\kappa_{i\ell}{\chi_\ell}_{|U_{ij}}-\sum_{\ell}\kappa_{j\ell}{\chi_\ell}_{|U_{ij}}\\
&=&\sum_{\ell}(\kappa_{i\ell}+\kappa_{\ell j}+\kappa_{ji}+\kappa_{ij}){\chi_\ell}_{|U_{ij}}\\
&=&\sum_{\ell}(\delta\tilde{\kappa}_{i\ell j}+\kappa_{ij}){\chi_\ell}_{|U_{ij}}\\
&=&\kappa_{ij},
\end{eqnarray*}
since $\delta\tilde{\kappa}=0$.

Its exterior differential $\tilde{\omega}:=d\tilde{\phi}$ satisfies $\delta\tilde{\omega}=\delta d\tilde{\phi}=d\delta\tilde{\phi}=d\tilde{\kappa}=0$, therefore it defines a global closed 1-form $\omega$. If $\kappa\in\ell^{p_1}$, $\omega\in L^{p_1}$ as well. If $\kappa=du$ where $u$ is an $\ell^q$ 0-cochain, the corresponding 0-form valued 0-cochain $\tilde{\phi}$ satisfies
\begin{eqnarray*}
\phi_i =\sum_{\ell}(u_{i}-u_{\ell})\chi_\ell=u_i-\psi
\end{eqnarray*}
where
\begin{eqnarray*}
\psi:=\epsilon(u)=\sum_{\ell}u_{\ell}\chi_\ell
\end{eqnarray*}
is a global $L^q$ 0-form. Furthermore, $\tilde{\omega}=d\tilde{\phi}=-d\psi$, so $\omega=-d\psi$. Thus we have a well-defined bounded linear map $L^{\pi}H^1(T)\to L^{\pi}H^1(X)$. 

The maps just defined in cohomology are inverses of each other. Indeed, all we have used is $\epsilon$ which inverts $\delta$ and Poincar\'e inequalities which allow to invert $d$. The map in one direction is $\delta\circ d^{-1}$; in the opposite direction, it is $d\circ \epsilon=d\circ\delta^{-1}$.

To sum up, the argument uses spaces of differential forms on open sets $U_i$'s and intersections $U_{ij}$'s, the operators $d$ and $\delta$, the inverse $\epsilon$ of $\delta$ provided by a partition of unity, the possibility to invert $d$ locally. The procedure $\delta\circ d^{-1}$ amounts to finding $b$ that relates a globally defined closed 1-form $\omega$ and a scalar 1-cocycle $\kappa$ via
\begin{eqnarray*}
(d+\delta)b=\tilde{\omega}+\tilde{\kappa},
\end{eqnarray*}
revealing the role played by the complex $\pm d+\delta$. Incidentally, we see that the inclusion $\ell^p\subset\ell^q$ when $p\leq q$ plays a role.

\subsection{The general case}

A bit of notation will help. Let $T$ be a simplicial complex, with a Banach space $\mathcal{B}_{i_0 \ldots i_h}$ attached at each simplex. Let $C^{h}(T,\mathcal{B})$ denote the set of cochains, i.e. skew-symmetric functions $\kappa$ on oriented simplices with values in $\mathcal{B}$ (i.e. $\kappa(i_0 \ldots i_h)$ is a vector in $\mathcal{B}_{i_0 \ldots i_h}$ for each simplex $i_0 \ldots i_h$).
Denote by
\begin{eqnarray*}
\mathcal{K}^{\pi,h}=\{\kappa\in C^{h}(T,\mathcal{B})\,;\,|\kappa|\in\ell^{p_h}\textrm{ and }|\delta\kappa|\in\ell^{p_{h+1}}\}.
\end{eqnarray*}

Let $C^{h,k}$ denote the space of $h$-cochains with values in $k$-forms, equipped with the $\mathcal{K}^{p_{h+k},h}(\Omega^{\pi,k})$-norm (here, $\mathcal{B}_{i_0 \ldots i_h}=\Omega^{\pi,k}(U_{i_0 \ldots i_h})$). It has two bounded differentials, $d'=\delta$ and $d''=(-1)^h d$, which anti-commute, thus $d'+d''$ is again a complex. Note that the space of globally defined, $L^{p_k}$ $k$-forms, is $\Omega^{\pi,k}(X)=C^{0,k}\cap\mathrm{ker}(d')$ and that the space of $\ell^{p_h}$ scalar valued $h$-cochains coincides with $C^{h,0}\cap\mathrm{ker}(d'')$. The choice of exponent $p_{h+k}$ in the definition of $C^{h,k}$, constant along diagonals of the bi-complex, makes it possible to iterate $\delta\circ d^{-1}$ and $d\circ\epsilon$ and relate $L^{p_k}$ and $\ell^{p_k}$ cohomologies.

Say a complex of Banach spaces $(C^{\cdot},d)$ is \emph{acyclic} up to degree $L$ if its cohomology vanishes in all degrees up to $L$.

We show that lines and columns of our bi-complex are acyclic.

\begin{lem}
If $p\leq q$, then $\|\cdot\|_{\ell^q}\leq \|\cdot\|_{\ell^p}$.
\end{lem}

\begin{pf}
If $x_i\geq 0$ and $\sum x_i=1$, then all $x_i$ are $\leq 1$, whence $\sum x_i^{q/p}\leq 1$. Applying this to
\begin{eqnarray*}
x_i=\frac{|a_i|^p}{\sum|a_j|^p}
\end{eqnarray*}
yields
\begin{eqnarray*}
\sum (|a_i|^p)^{q/p}\leq (\sum|a_j|^p)^{q/p},
\end{eqnarray*}
whence the announced inequality.
\end{pf}

\begin{lem}
\label{pq}
Let $\mathcal{U}=(U_i)_i$ be a open covering of a Riemannian manifold. Assume that the volumes of $U_i$'s are bounded, and that $\mathcal{U}$ admits a partition of unity $(\chi_i)_i$ with uniformly bounded Lipschitz constants. Given $\phi\in C^{h,k}$, i.e. $\phi$ is the data of a differential $k$-form on each $U_{i_0 \ldots i_h}$, set
\begin{eqnarray*}
\epsilon(\phi)_{i_0 \ldots i_{h-1}}=\sum_{j}\chi_j \phi_{j i_0 \ldots i_{h-1}}.
\end{eqnarray*}
Then $\epsilon:C^{h,k}\to C^{h-1,k}$ is bounded and $1=\epsilon\delta+\delta\epsilon$.
\end{lem}

\begin{pf}
Computing
\begin{eqnarray*}
((\epsilon\delta+\delta\epsilon)\phi)_{i_0 \ldots i_h}
&=&\sum_{j}\chi_j(\delta\phi)_{j i_0 \ldots i_h}+\sum_{\ell}(-1)^{\ell}(\epsilon\phi)_{i_0 \ldots \widehat{i_\ell}\ldots i_h}\\
&=&\sum_{j}\chi_j(\phi_{i_0 \ldots i_h}+\sum_{\ell}(-1)^{\ell+1}\phi_{j i_0 \ldots \widehat{i_\ell}\ldots i_h})\\
&&+\sum_{\ell}(-1)^{\ell}\sum_j \chi_j \phi_{j i_0 \ldots \widehat{i_\ell}\ldots i_h}\\
&=&\sum_{j}\chi_j\phi_{i_0 \ldots i_h}\\
&=&\phi_{i_0 \ldots i_h},
\end{eqnarray*}
shows that $\epsilon\delta+\delta\epsilon=1$.

Multiplying a differential $k$-form $\phi_{i_0 \ldots i_{h-1}j}$ with $\chi_j$ does not increase its $L^{p_k}$-norm. For its differential,
\begin{eqnarray*}
\|d(\chi_j \phi_{i_0 \ldots i_{h-1}j})\|_{L^{p_{k+1}}(U_{i_0 \ldots i_{h-1}j})}&\leq& \|d\phi_{i_0 \ldots i_{h-1}j}\|_{L^{p_{k+1}}(U_{i_0 \ldots i_{h-1}j})}\\
&&+L\|\phi_{i_0 \ldots i_{h-1}j}\|_{L^{p_{k+1}}(U_{i_0 \ldots i_{h-1}j})}.
\end{eqnarray*}
where $L$ is a Lipschitz bound on $\chi_j$. By H\"older's inequality, since $p_{k+1}\leq p_k$, the second term is bounded above by $\|\phi_{i_0 \ldots i_{h-1}j}\|_{L^{p_{k}}(U_{i_0 \ldots i_{h-1}j})}$ times a power of the volume of $U_{j}$, which is assumed to be bounded above. Thus multiplication with $\chi_j$ is continuous in local $\Omega^{\pi,k}$ norms. 

Since $\chi_j$ has compact support in $U_j$, $\chi_j \phi_{i_0 \ldots i_{h-1}j}$ extends by 0 to $U_{i_0 \ldots i_{h-1}}$ without increasing its $\Omega^{\pi,k}$ norm. Therefore $\epsilon$ is bounded from $\ell^{p_h}(\Omega^{\pi,k})$ to $\ell^{p_h}(\Omega^{\pi,k})$, and thus from $\ell^{p_h}(\Omega^{\pi,k})$ to $\ell^{p_{h-1}}(\Omega^{\pi,k})$ by Lemma \ref{pq}. With the identity $1=\epsilon\delta+\delta\epsilon$, this shows that $\epsilon:C^{h,k}\to C^{h-1,k}$ is bounded. 
\end{pf}

\begin{cor}
Under the asumptions of Lemma \ref{pq}, the horizontal complexes $(C^{\cdot,k},\delta)$ are acyclic.
\end{cor}

\begin{lem}
\label{poin}
Assume that $\pi$ satisfies $1<p_h<\infty$ for all $h$. Assume that the open covering $\mathcal{U}$ satisfies the following uniformity property, for some constant $M$: for each nonempty intersection $U_{i_0\ldots i_h}$, there is a $M$-bi-Lipschitz homeomorphism of $U_{i_0\ldots i_h}$ to the unit ball of $\R^n$. Then the vertical complexes $(C^{h,\cdot},d)$ are acyclic.
\end{lem}

\begin{pf}
We use the fact that inequality $(Poincare_{q,p}(k))$ is valid for $\lambda=1$ (no loss on the size of domain) if $p>1$ and $q<\infty$. This is due to Iwaniec and Lutoborsky, \cite{IL}: Cartan's formula provides an explicit inverse to $d$ on bounded convex domains, which is bounded from $L^p$, $p>1$, to $L^q$, $q<\infty$, provided $\frac{1}{p}-\frac{1}{q}\leq\frac{1}{n}$.
\end{pf}

\begin{lem}
\label{leray}
Let $(C^{h,k},d',d'')$ be a bi-complex of Banach spaces indexed by $\mathbb{N}\times\mathbb{N}$. Assume that the horizontal complexes $(C^{\cdot,k},d')$ all are acyclic up to degree $L$. Then the inclusion of $(C^{0,\cdot}\cap \mathrm{ker}(d'),d'')$ into $(C^{\cdot,\cdot},d'+d'')$ induces an isomorphism in cohomology up to degree $L$.
\end{lem}

\begin{pf}
Let us replace $C^{h,k}$ with $C^{h,k}\cap\mathrm{ker}(d')$ when $h+k=L$ and by 0 when $h+k>L$. This does not affect the conclusion in degree $L$, and allows to assume acyclicity in all degrees.

If $a\in C^{\cdot,\cdot}$, let $a_m$ denote the sum of the components of $a$ of $h$-degree equal to $m$. For an integer $\ell$, let
\begin{eqnarray*}
A^\ell&:=&\{a\in C^{\cdot\leq \ell,\cdot}\,;\,(d'+d'')a\in C^{\cdot\leq \ell,\cdot}\}\\
&=&\{a\in C^{\cdot\leq \ell,\cdot}\,;\,d'a_\ell=0\}.
\end{eqnarray*}
Then $(A^\ell,d'+d'')$ is a subcomplex. 

One shows that for all $\ell$, the inclusion of $A^\ell$ into $A^{\ell+1}$ induces an isomorphism in cohomology. If $a\in A^{\ell+1}$ is $d'+d''$-closed, by acyclicity, there exists $b\in C^{\cdot\leq \ell,\cdot}$ such that $d'b=a_{\ell+1}$. Then $a'=a-(d'+d'')b\in C^{\cdot\leq \ell,\cdot}$, it is $d'+d''$-closed, so $a'\in A^\ell$. This shows that the inclusion $A^\ell\to A^{\ell+1}$ is onto in cohomology. Let $a\in A^\ell\cap\mathrm{ker}(d'+d'')$. Assume that there exists $b\in A^{\ell+1}$ such that $a=(d'+d'')b$. By acyclicity, there exists $c\in C^{\cdot\leq \ell,\cdot}$ such that $d'c=b_{\ell+1}$. Then $b':=b-(d'+d'')c\in C^{\cdot\leq \ell,\cdot}$. Since
\begin{eqnarray*}
(d'+d'')b'=(d'+d'')b=a\in C^{\cdot\leq \ell,\cdot},
\end{eqnarray*}
$b'\in A^\ell$. Since $a=(d'+d'')b'$, this shows that the inclusion $A^\ell\to A^{\ell+1}$ is injective in cohomology.
\end{pf}

\begin{cor}
Under the asumptions of Lemmas \ref{pq} and \ref{poin}, $L^\pi$-cohomology of differential forms and $\ell^\pi$-cohomology of simplicial cochains of the nerve coincide.
\end{cor}

\begin{pf}
Apply Lemma \ref{leray} twice, once with $d'=\delta$ and $d''=\pm d$, once with $d'=d$, $d''=\pm\delta$.
\end{pf}

\begin{rem}
Assume an even stronger form of Poincar\'e inequality holds: up to degree $L$, there exist bounded linear operators $e:\Omega^{\pi,k}(U_{i_0\ldots i_h})\to \Omega^{\pi,k-1}(U_{i_0\ldots i_h})$ with uniformly bounded norms such that $1=de+ed$. Then the conclusion is stronger: there exists a homotopy of complexes $\Omega^{\pi,\cdot}(X)\to\mathcal{K}^{\pi,\cdot}(T)$ up to degree $L$.
\end{rem}

\section{A customized version of Leray's acyclic coverings theorem}
\label{leray cust}

The above argument has the following drawbacks:
\begin{itemize}
  \item It requires Poincar\'e inequalities without loss on the size of domain, which are not known in all cases.
  \item It makes strong assumptions on coverings, see Lemma \ref{poin}.
\end{itemize}
Fortunately, a modification of the homological algebra allows for weaker assumptions on coverings and for weaker Poincar\'e inequalities, allowing loss on the size of domain, as stated in Theorem \ref{poincare}. 

\subsection{Existence of uniform coverings}

\begin{dfi}
Let $X$ be a metric space. A \emph{uniform sequence of nested coverings} is a sequence $\mathcal{U}^{0},\ldots,\mathcal{U}^{L}$ of open coverings of $X$ with the following properties, for some constants $M>0$, $r>0$ and some model pair of metric spaces $Z'\subset Z$,
\begin{enumerate}
  \item Nesting: for each $i$ and all $j=1,\cdots,L$, $U_i^{j-1}\subset U_i^{j}$.
  \item Bounded size: the diameters of $U_i^L$'s are bounded; each $U_i^0$ contains a ball of radius $r$, and these balls are disjoint.
  \item Bounded multiplicity: every point of $X$ is contained in at most $M$ open sets $U_i^L$.
  \item Bounded partition of unity: $\mathcal{U}^0$ has a partition of unity with bounded Lipschitz constants.
  \item Contractibility: each $U_i^{0}$ is contractible within $U_i^{1}$.
  \item Model: for each pair $(U_{i_0\ldots i_h}^{j-1},U_{i_0\ldots i_h}^{j})$ such that $U_{i_0\ldots i_h}^0$ is nonempty, there is an $M$-bi-Lipschitz map $\phi_{i_0\ldots i_h,j}:Z\to X$ such that $\phi_{i_0\ldots i_h,j}(Z)\subset U_{i_0\ldots i_h}^{j}$ and $U_{i_0\ldots i_h}^{j-1}\subset\phi_{i_0\ldots i_h,j}(Z')$.
\end{enumerate}
\end{dfi}

With some lead over Section \ref{contact}, let us define bounded geometry in the contact sub-Riemannian case.
\begin{dfi}
\label{defbddgeocontact}
Let $X$ be a contact manifold equipped with a sub-Rieman\-nian metric. Say that $X$ has bounded geometry if there exist $M>0$ and $r_0>0$ and for every point $x$ a smooth contactomorphism $\phi_x$ of the unit ball of $\H^m$ to an open subset of $X$, mapping the origin to $x$, and such that $B(x,r_0)\subset\phi_x(B(1))$, and $\phi_x$ is $M$-bi-Lipschitz.
\end{dfi}

\begin{pro}
\label{existcover}
Let $X$ be a bounded geometry Riemannian or contact manifold. Then $X$ admits uniform sequences of nested coverings of arbitrary length, where the models are pairs of concentric Euclidean (resp. Heisenberg) balls whose ratio of radii can be chosen arbitrarily.
\end{pro}

\begin{pf}
Fix $\lambda>1$. Let $r_0$ be the radius occurring in the definition of bounded geometry. Let $r=(\lambda M)^{-2L}r_0$. Pick a maximal packing of $X$ by disjoint $r$-balls. Let $\mathcal{U}^0$ be the collection of twice larger balls $U_i^0=B(x_i,2r)$, and $U_i^j=B(x_i,2(\lambda M)^{2j} r)$. The nesting, size and multiplicity requirements are satisfied. The partition of unity can be constructed from the distance function to $x_i$, it is uniformly Lipschitz.

If $U_{i_0\ldots i_h}^0$ is nonempty, then all $x_{i_\ell}$ belong to $B(x_{i_0},4r)$, thus 
\begin{eqnarray*}
B(x_{i_0},((\lambda M)^{2j}-4)r)\subset U_{i_0\ldots i_h}^j.
\end{eqnarray*}
Consider given chart $\phi_{x_{i_0}}:B(1)\to X$, whose image contains $U_{i_0\ldots i_h}^L$ by construction. Then 
\begin{eqnarray*}
\phi_{x_{i_0}}(B(\frac{1}{M}((\lambda M)^{2j}-4)r))\subset B(x_{i_0},((\lambda M)^{2j}-4)r)\subset U_{i_0\ldots i_h}^j.
\end{eqnarray*}
On the other hand,
\begin{eqnarray*}
\phi_{x_{i_0}}^{-1}(B(x_{i_0},(\lambda M)^{2j-2}r))\subset B(M(\lambda M)^{2j-2}r),
\end{eqnarray*}
thus
\begin{eqnarray*}
U_{i_0\ldots i_h}^{j-1}\subset B(x_{i_0},(\lambda M)^{2j-2}r)\subset \phi_{x_{i_0}}(B(M(\lambda M)^{2j-2}r)).
\end{eqnarray*}
Thus one can set $\phi_{i_0\ldots i_h,j}=\phi_{x_{i_0}}$ precomposed with dilation (in Euclidean space or Heisenberg group) by factor $M(\lambda M)^{2j-2}$. The pair $Z'=B(r)$, $Z=B(\lambda r)$ of concentric balls serves as a model. Indeed, the ratio of radii
\begin{eqnarray*}
\frac{\frac{1}{M}(\lambda M)^{2j}-4}{M(\lambda M)^{2j-2}}=\lambda^2-\frac{4}{\lambda^{2j-2}M^{2j-1}}\geq\lambda^2 -\frac{4}{M}\geq \lambda,
\end{eqnarray*}
provided $M\geq 2$ and $\lambda\geq 2$. 

The contractibility requirement follows from the existence of model, since model balls are contractible.
\end{pf}

\begin{rem}
By definition, in a uniform sequence of nested coverings, Poincar\'e inequalities as in Theorem \ref{poincare} hold with uniform constants for all pairs $(U_{i_0\ldots i_h}^{j-1},U_{i_0\ldots i_h}^{j})$ such that $U_{i_0\ldots i_h}^0$ is nonempty.
\end{rem}
Indeed, pull-back by $M$-bi-Lipschitz diffeomorphisms (resp. contactomorphisms) expands or contracts $L^p$ norms of differential forms (resp. Rumin forms) by at most a power of $M$. This is also true for $BMO$ and $\mathcal{H}^1$, \cite{BN}.

\subsection{Closed 1-forms and 1-cocycles}

Let $\mathcal{U}^0$ and $\mathcal{U}^1$ be nested open coverings, i.e. for all $i$, $U_i^0 \subset U_i^1$. One assumes that Poincar\'e inequality applies to each pair $(U_i^0,U_i^1)$. One introduces the two bi-complexes $C^{\cdot,\cdot,0}$ and $C^{\cdot,\cdot,1}$ associated with the two coverings. The simplicial complexes $T^0$ and $T^1$ share the same vertices, but $T^0$ has less simplices. Without change in notation, let us associate the trivial vectorspace to simplices of $T^1$ which do not belong to $T^0$. Let $r:C^{\cdot,\cdot,1}\to C^{\cdot,\cdot,0}$ denote the restriction operator (which vanishes for empty intersections $U^0_{i_0\ldots i_h}$). It commutes with $d$ and $\delta$.

Using covering $\mathcal{U}^1$, a globally defined closed 1-form $\omega$ on $X$ defines an element $\tilde{\omega}^1\in C^{0,1,1}$. The primitive $\tilde{\phi}^0\in C^{0,0,0}$ provided by local Poincar\'e inequalities satisfies $d\tilde{\phi}^0=r(\tilde{\omega}^1)$. $\tilde{\kappa}^0=\delta(\tilde{\phi}^0)$ defines a 1-cocycle of $\mathcal{U}^0$. The inverse procedure, from 1-cocycles of $\mathcal{U}^0$ to closed 1-forms, is unaffected by covering $\mathcal{U}^1$. Both procedures, when precomposed with the restriction operator $r$, coincide with the procedures defined earlier, i.e. provide the required cohomology isomorphism relative to covering $\mathcal{U}^0$.

To sum up, only one simplicial complex is needed, the nerve of the covering $\mathcal{U}^0$ by small open sets.

\subsection{General case}

One starts with a uniform sequence $(\mathcal{U}^{j})_{j=0,\ldots,L}$ of nested coverings. Let $C^{\cdot,\cdot,j}$ denote the bi-complex of cochains of $T^0$ with values in differential forms on intersections $U_{i_0 \ldots i_h}^j$ (this differs from the bi-complex associated to $\mathcal{U}^j$). There are restriction maps $r:C^{\cdot,\cdot,j}\to C^{\cdot,\cdot,j-1}$ which commute with $d$ and $\delta$. Poincar\'e inequalities state that a form of acyclicity holds: $r$ induces the 0 map in cohomology.

\begin{dfi}
Let $(r:C^{\cdot,j}\to C^{\cdot,j-1},d)$ be a commutative diagram of complexes. Say the diagram is \emph{$r$-acyclic} up to degree $L$ if $r$ induces the 0 map in cohomology up to degree $L$.
\end{dfi}

Lemma \ref{leray} is replaced with

\begin{lem}
\label{leraybis}
Let $(r:C^{\cdot,\cdot,j}\to C^{\cdot,\cdot,j-1},d',d'')$ be a commutative diagram of  bi-complexes. Assume that all horizontal diagrams $(r:C^{\cdot,k,j}\to C^{\cdot,k,j-1},d')$, are $r$-acyclic up to degree $L$. Let $\iota_j$ denote the cohomology map induced by the inclusion of $(C^{0,\cdot,j}\cap \mathrm{ker}(d''),d')$ into $(C^{\cdot,\cdot,j},d'+d'')$. Then up to degree $L$, the image of $\iota_0$ contains the image of $r^L$, and the kernel of $\iota_L$ is contained in the kernel of $r^L$.
\end{lem}

\begin{pf}
As before, one may assume that the bi-complex has finitely many terms and is $r$-acyclic in all degrees. Denote by $A^{\ell,j}$ the sub-complexes introduced in the proof of Lemma \ref{leray}, relative to the $j$-th complex, i.e.
\begin{eqnarray*}
A^{\ell,j}&:=&\{a\in C^{\cdot\leq \ell,\cdot, j}\,;\,(d'+d'')a\in C^{\cdot\leq \ell,\cdot,j}\}.
\end{eqnarray*}
The same argument as in the proof of Lemma \ref{leray} shows that, for all $\ell$,
\begin{enumerate}
  \item for all closed $a\in A^{\ell+1,j}$, there exists $a'\in A^{\ell,j-1}$ such that $ra-a'\in(d'+d'')(A^{\ell+1,j})$.
  \item if $a\in A^{\ell,j}$ belongs to $(d'+d'')(A^{\ell+1,j})$, then $ra\in(d'+d'')(A^{\ell,j-1})$.
\end{enumerate}
It suffices to iterate $L$ times to obtain the claimed statement. Indeed, each time $d'$ is inverted, degree decreases by 1, so at most $L$ inversions are required. 
\end{pf}

\begin{cor}
\label{deRham}
Let $X$ be a bounded geometry Riemannian manifold. Pick a uniform sequence of nested coverings of length $2L$. Assume that $L^\pi$-cohomology is modified as prescribed in Theorem \ref{qib} for exceptional values of $(p,q,k)$. Then $L^\pi$-cohomology of differential forms and $\ell^\pi$-cohomology of simplicial cochains of the nerve of $\mathcal{U}^0$ are isomorphic as vectorspaces. The isomorphism maps the exact cohomology of $X$ to the exact cohomology of the nerve.
\end{cor}

\begin{pf}
The following diagram commutes.
\begin{eqnarray*}
\begin{CD}
H^{\cdot}(C^{0,\cdot,0}\cap \mathrm{ker}(d'')) @>{\iota_0}>> H^{\cdot}(C^{\cdot,\cdot,0}) @<{\iota'_0}<< H^{\cdot}(C^{\cdot,0,0}\cap \mathrm{ker}(d'))\\
@A{\rho^L}A\simeq A @A{R^L}AA @A{\rho'^L}A\simeq A \\
H^{\cdot}(C^{0,\cdot,L}\cap \mathrm{ker}(d'')) @>{\iota_L}>> H^{\cdot}(C^{\cdot,\cdot,L}) @<{\iota'_L}<< H^{\cdot}(C^{\cdot,0,L}\cap \mathrm{ker}(d'))\\
@A{\rho^{L}}A\simeq A @A{R^{L}}AA @A{\rho'^{L}}A\simeq A \\
H^{\cdot}(C^{0,\cdot,2L}\cap \mathrm{ker}(d'')) @>{\iota_{2L}}>> H^{\cdot}(C^{\cdot,\cdot,2L}) @<{\iota'_{2L}}<< H^{\cdot}(C^{\cdot,0,2L}\cap \mathrm{ker}(d'))\\\end{CD}
\end{eqnarray*}
For clarity, we used different notations, $\rho$, $R$ and $\rho'$, for the cohomology maps induced by $r$ for the 3 different complexes.

1. The $d''=\delta$ complexes are $r$-acyclic (in fact, acyclic in the usual sense, but we do not need this favourable circumstance). The complex $C^{\cdot,0,j}\cap \mathrm{ker}(d')$ consists of scalar simplicial cochains of the nerve $T^0$, restriction $r$ has no effect on them. Therefore the cohomology map $\rho'$ between consecutive levels is an isomorphism. From Lemma \ref{leraybis}, it follows that the image $I'_0$ of $\iota'_0:H^{\cdot}(C^{\cdot,0,0}\cap \mathrm{ker}(d'),d'')\to H^{\cdot}(C^{\cdot,\cdot,0},d'+d'')$ composed with ${\rho'}^L$ contains the image $I$ of $R^L:H^{\cdot}(C^{\cdot,\cdot,L},d'+d'')\to H^{\cdot}(C^{\cdot,\cdot,0},d'+d'')$. Also, $\iota'_L:H^{\cdot}(C^{\cdot,0,L}\cap \mathrm{ker}(d'),d'')\to H^{\cdot}(C^{\cdot,\cdot,L},d'+d'')$ is injective. Let $I'_L\subset H^{\cdot}(C^{\cdot,\cdot,L},d'+d'')$ denote its image.

2. Thanks to Theorem \ref{poincare}, the $d'=\pm d$ complexes are $r$-acyclic. The complex $C^{0,\cdot,j}\cap \mathrm{ker}(d'')$ consists of globally defined differential forms, restriction $r$ has no effect on them. Therefore the cohomology map $\rho$ between consecutive levels is an isomorphism. Lemma \ref{leraybis} implies that the image $I_0$ of $\iota_0:H^{\cdot}(C^{0,\cdot,0}\cap \mathrm{ker}(d''),d')\to H^{\cdot}(C^{\cdot,\cdot,0},d'+d'')$ composed with ${\rho}^L$ contains $I$, and that $\iota_L:H^{\cdot}(C^{0,\cdot,L}\cap \mathrm{ker}(d''),d')\to H^{\cdot}(C^{\cdot,\cdot,L},d'+d'')$ is injective. Let $I_L\subset H^{\cdot}(C^{\cdot,\cdot,L},d'+d'')$ denote its image.

3. Since $R^L\circ \iota_L\circ (\rho^{-1})^L=\iota_0$, $I_0=\mathrm{im}(\iota_0)\subset \mathrm{im}(R^L)=I$. We conclude that $I_0=I$. Similarly, $I=I'_0$, hence $I_0=I'_0$. For the same reason, using the bottom part of the diagram, $I_L=I'_L$.

Therefore $(\iota'_L)^{-1}\circ\iota_L$ is a bijection 
\begin{eqnarray*}
L^\pi H^\cdot(X)=H^{\cdot}(C^{0,\cdot,L}\cap \mathrm{ker}(d''))\to H^{\cdot}(C^{\cdot,0,L}\cap \mathrm{ker}(d'))=\ell^\pi H^\cdot(T^0).
\end{eqnarray*}

4. The construction provides an isomorphism between quotients of spaces of forms/cochains of finite norms, but also between quotients of larger spaces of forms/cochains without any decay condition. Therefore the isomorphism is compatible with forgetful maps $H^{q,p,k}\to H^k$ to ordinary (un-normed) cohomology, it maps kernel to kernel, exact cohomology to exact cohomology.
\end{pf}

\begin{rem}
Since $d''=\delta$ has a bounded linear inverse $\epsilon$, i.e. $1=\delta\epsilon+\epsilon\delta$, the map $\iota'_L$ has a continuous inverse, hence the linear isomorphism $(\iota'_L)^{-1}\circ\iota_L$ is continuous. If $\ell^\pi H^\cdot(T^0)$ is Hausdorff, so is $L^\pi H^\cdot(X)$ and both are isomorphic as Banach spaces.
\end{rem}

\begin{rem}
\label{specialcases}
Assume a slightly stronger form of Poincar\'e inequality holds: up to degree $L$, there exist bounded linear operators 
$$
e:\Omega^{\pi,k}(U^{j}_{i_0\ldots i_h})\to \Omega^{\pi,k-1}(U^{j-1}_{i_0\ldots i_h})
$$
with uniformly bounded norms such that $1=de+ed$. Then the conclusion is stronger: there exists a homotopy of complexes $\Omega^{\pi,\cdot}(X)\to\mathcal{K}^{\pi,\cdot}(T)$ up to degree $L$. In particular, the vectorspace isomorphism is topological, it induces isomorphisms of reduced cohomology and torsion. 
\end{rem}
The stronger assumption holds unless $p=1$ and $q=\frac{n}{n-1}$, or $p=n$ and $q=\infty$, \cite{IL}. It fails if $p=n$, $q=\infty$ (\cite{BB1}, Proposition 2, for $k=n$, \cite{BB2}, Proposition 9, for other values of $k$).

\subsection{Limiting cases}

To show that $L^{q,p} H^k(X)$ and its discretized version $\ell^{q,p} H^k(T^0)$ are isomorphic, one defines a vector $\pi$ by $p_0=\cdots=p_{k-1}=q$ and $p_k=p$. If a limiting case arises, i.e. $\frac{1}{p}-\frac{1}{q}=\frac{1}{n}$ and either $p=1$ or $q=\infty$, it is only in degree $k$ that a limiting Poincar\'e inequality is required. This is why the restrictions on $(p,q,k)$ appearing in Theorem \ref{poincare} are exactly reflected in Theorem \ref{qib}.
 
\subsection{Multiplicative structure}

Differential forms form a graded differential algebra: the wedge product satisfies
\begin{eqnarray*}
d(\alpha\wedge\beta)=d\alpha\wedge\beta+(-1)^{\mathrm{deg}(\alpha)}\alpha\wedge d\beta.
\end{eqnarray*}
Simplicial cochains do as well: the cup-product satisfies
\begin{eqnarray*}
\delta(\kappa\smile\lambda)=\delta\kappa\smile\lambda+(-1)^{\mathrm{deg}(\kappa)}\kappa\smile \delta\lambda.
\end{eqnarray*}
The tensor product of two graded differential algebras inherits the structure of a graded differential algebra (the algebra of differential forms on the product of two manifolds illustrates this). Therefore, if $\phi\in C^{h,k}$ and $\phi'\in C^{h',k'}$, set
\begin{eqnarray*}
(\phi\smile\phi')_{i_0\ldots i_{h+h'}}=&\\
(-1)^{kh'}\sum_{\sigma\in\mathfrak{S}_{h,h'}}(-1)^\sigma &(\phi_{\sigma(i_0)\ldots \sigma(i_h)}\wedge\phi'_{\sigma(i_h)\ldots \sigma(i_{h+h'})})_{|U_{i_0 \ldots i_{h+h'}}},
\end{eqnarray*}
where $\mathfrak{S}_{h,h'}$ denotes the set of permutations of $\{0,\ldots,h+h'\}$ which are increasing on $\{0,\ldots,h\}$ and on $\{h,\ldots,h+h'\}$. Set
\begin{eqnarray*}
(d'+d'')\phi=\delta\phi+(-1)^h d\phi.
\end{eqnarray*}
The multiplication is continuous $L^p\otimes L^{p'}\to L^{p''}$ provided
$\frac{1}{p}+\frac{1}{p'}=\frac{1}{p''}$. It maps $L^q \otimes L^{p'}\to L^{q''}$ and $L^{q'} \otimes L^{p}\to L^{q''}$ provided $\frac{1}{p}+\frac{1}{q'}=\frac{1}{q''}$ and $\frac{1}{p'}+\frac{1}{q}=\frac{1}{q''}$.

This multiplication descends to cohomology and restricts to the usual cup-product on de Rham and simplicial cohomology. Since the isomorphisms of  Lemmas \ref{leray} and \ref{leraybis} arise from multiplication preserving inclusions, they preserve multiplication.

One concludes that, provided
\begin{eqnarray*}
\frac{1}{p}+\frac{1}{p'}=\frac{1}{p''}\quad\textrm{and}\quad \frac{1}{p}+\frac{1}{q'}=\frac{1}{q''}=\frac{1}{p'}+\frac{1}{q},
\end{eqnarray*}
the cup-product
\begin{eqnarray*}
\smile:H^{q,p,k}(X)\otimes H^{q',p',k'}(X)\to H^{q'',p'',k+k'}(X)
\end{eqnarray*}
is well-defined, and can be computed either using differential forms or simplicial cochains.

From now on, we shall work with the simplicial complex $T$ and its 0-skeleton $Y$.

\section{Leray's theorem for simplicial complexes}
\label{leray simpl}

Next, we establish an analogue of Corollary \ref{deRham} where manifolds are replaced with simplicial complexes and differential forms with simplicial cochains. The analytic point, Poincar\'e inequalities for differential forms, turns out to be replaced with a purely topological fact.

In this section, radii $R$ are integers. All simplicial complexes have bounded geometry, i.e. every vertex belongs to a bounded number of simplices. The exponent sequence $p_0 \geq \cdots \geq p_k \geq \cdots$ is nonincreasing.

\subsection{Uniform vanishing of cohomology}

\begin{dfi}
Say a simplicial complex $X$ with 0-skeleton $Y$ has uniformly vanishing cohomology up to degree $L$ if for every $R>0$, there exists $\tilde{R}(R)$ such that for every point $y\in Y$, the maps $H^{j}(B(y,\tilde{R}))\to H^{j}(B(y,R))$ induced by inclusion $B(y,R)\subset B(y,\tilde{R})$ vanish for all $j=0,\ldots,L$.
\end{dfi}

\begin{exa}
Assumption holds if $X$ has vanishing cohomology up to degree $L$ and a cocompact automorphism group.
\end{exa}
Indeed, by duality, the assumption is that homology vanishes. The vectorspace of cycles in $B(y,R)$ is finite dimensional. Pick a finite basis. Every element bounds a finite chain, all these chains are contained in some ball $B(y,\tilde{R})$. Thus all maps $H^{j}(B(y,\tilde{R}))\to H^{j}(B(y,R))$ vanish. $\tilde{R}$ depends on $R$ and $y$. If $X$ has a cocompact automorphism group, $\tilde{R}$ depends on $R$ only.

\subsection{Poincar\'e inequality for simplicial complexes}

\begin{lem}
\label{SPoincare}
Let $X$ be a simplicial complex $X$ with uniformly vanishing cohomology up to degree $L$. Then Poincar\'e inequalities hold for all pairs $(B(y,R-1),B(y,\tilde{R}(R)))$ up to degree $L$. For the subspace of exact cocycles, Poincar\'e inequalities hold for all degrees. In both cases, constants do not depend on $y$
\end{lem}

\begin{pf}
Let $C(y,R)$ (resp. $C'(y,R)$) be the union of simplices contained in $B(y,R)$ (resp. intersecting $B(y,\tilde{R})$). As $y$ varies, at most finitely many different pairs of complexes $(C,C')$ are encountered. By assumption, for each pair, the cohomology maps $H^j(C')\to H^j(C)$ vanish if $j\leq L$. If $j>L$, the cohomology maps $EH^j(C')\to EH^j(C)$ vanish by definition. Since simplicial cochains of $C$ and $C'$ form finite dimensional vectorspaces, Poincar\'e inequality is nothing more than this vanishing. Uniformity of constants arises from finiteness of the collection of maps.
\end{pf}

\subsection{Vertical $r$-acyclicity}

Using uniform vanishing of cohomology, one constructs nested coverings as follows. Fix $R\in\mathbb{N}$. The specifications are that all $U_i^0$ be $R$-balls and each pair $(U_{i_0 \ldots i_h}^{j},U_{i_0 \ldots i_h}^{j+1})$ such that $U_{i_0 \ldots i_h}^{0}\not=\emptyset$ satisfies Poincar\'e inequality.

Let $R_0=R$. Let covering $\mathcal{U}^0$ consist of subcomplexes $C(y,R_0)$, $y\in Y$. Pick $R_1=\tilde{R}(R_0+R)+R$ according to uniform cohomology vanishing, and let covering $\mathcal{U}^1$ consist of subcomplexes $C(y,R_1)$, $y\in Y$, set $R_2=\tilde{R}(R_1+R)+R$, and so on. If $y\in U_{i_0 \ldots i_h}^{0}$, then the centers of all $U_{i_{\ell}}^j$, $\ell=0,\ldots,h$, all $j$, belong to $B(y,R)$, so $U_{i_0 \ldots i_h}^{j}$ is contained in $B(y,R_j +R)$ and $U_{i_0 \ldots i_h}^{j+1}$ contains $B(y,R_{j+1}-R)$. Since $R_{j+1}-R\geq \tilde{R}(R_{j}+R)$, the pair $(B(y,R_{j}+R),B(y,R_{j+1}-R))$ satisfies Poincar\'e inequality. This shows that all relevant pairs $(U_{i_0 \ldots i_h}^{j},U_{i_0 \ldots i_h}^{j+1})$ satisfy Poincar\'e inequality. All other boundedness properties follow from the fact that $X$ has bounded geometry.

We consider the bi-complexes $C^{h,k,j}$ consisting of $h$-cochains of the nerve of $\mathcal{U}^0$ with values in $k$-cochains of intersections of open sets $U_{i_0 \ldots i_h}^j$ of $\mathcal{U}^j$. We truncate it: if $h+k>L+1$, we set $C^{h,k,j}=0$ and replace $\bigoplus_{h+k=L+1}C^{h,k,j}$ with its subspace of exact elements. Here, $d'=\delta$ is the covering coboundary, $d'=(-1)^h d$ is the simplicial coboundary of $X$. Let $r:C^{\cdot,\cdot,j}\to C^{\cdot,\cdot,j-1}$ denote the restriction map. From Lemma \ref{SPoincare}, vertical complexes are $r$-acyclic.

\subsection{Horizontal acyclicity}

\begin{lem}
The horizontal complexes $d':C^{\cdot,k,j}\to C^{\cdot+1,k,j}$ are acyclic.
\end{lem}

\begin{pf}
The same operator $\epsilon$ which inverts $\delta$ will be used for all coverings $\mathcal{U}^j$. It is made from a partition of unity $(\chi_i)$ for $\mathcal{U}^0$. Let $\eta_i$ denote the function on $Y$ which is 1 on $Y\cap U_i$ and 0 elsewhere. Set 
\begin{eqnarray*}
\chi_i=\frac{\eta_i}{\sum_{j}\eta_j}.
\end{eqnarray*}
If $\kappa$ a $k$-cochain on $T$, view $\chi_i$ as a 0-cochain and use the unskewsymmetrized cup-product to multiply $\chi_i$ with $\kappa$,
\begin{eqnarray*}
(\chi_i\smile\kappa)_{y_0 \ldots y_k}=\chi_i(y_0)\kappa_{y_0 \ldots y_k}.
\end{eqnarray*}
If $\kappa_{j i_0 \ldots i_{h-1}}$ is defined on $U_{j i_0 \ldots i_{h-1}}$, $\chi_j\smile\kappa$ extends by 0 to $U_{i_0 \ldots i_{h-1}}$. So the following $k$-cochain 
\begin{eqnarray*}
(\epsilon\kappa)_{i_0 \ldots i_{h-1}}=\sum_{j}\chi_j\smile\kappa_{j i_0 \ldots i_{h-1}}
\end{eqnarray*}
is well-defined on $U_{i_0 \ldots i_{h-1}}$. The identity $\epsilon\delta+\delta\epsilon=1$ is formal. The formula
\begin{eqnarray*}
d(\chi_j\smile\kappa_{j i_0 \ldots i_{h-1}})=(d\chi_j)\smile\kappa_{j i_0 \ldots i_{h-1}}+\chi_j\smile(d\kappa_{j i_0 \ldots i_{h-1}})
\end{eqnarray*}
shows that the norm of $(\epsilon\kappa)_{i_0 \ldots i_{h-1}}$ in $\ell^{p_k}(U_{i_0 \ldots i_{h-1}})$ is controlled by the norms of $\kappa_{j i_0 \ldots i_{h-1}}$ and $d\kappa_{j i_0 \ldots i_{h-1}}$ in $\ell^{p_k}(U_{j i_0 \ldots i_{h-1}})$. By H\"older's inequality, one can replace the latter by the norm of $d\kappa_{j i_0 \ldots i_{h-1}}$ in $\ell^{p_{k+1}}(U_{j i_0 \ldots i_{h-1}})$ (since the number of $k+1$-simplices in $U_{j i_0 \ldots i_{h-1}}$ is bounded). This shows that $\epsilon$ is bounded in local $\mathcal{K}^{\pi,k}$-norms. Adding terms up shows that $\epsilon$ is bounded from $\ell^{p_h}(\mathcal{K}^{\pi,k})$ to $\ell^{p_h}(\mathcal{K}^{\pi,k})$, and thus from $\ell^{p_h}(\mathcal{K}^{\pi,k})$ to $\ell^{p_{h-1}}(\mathcal{K}^{\pi,k})$ by Lemma \ref{pq}. With the identity $\epsilon\delta+\delta\epsilon=1$, this shows that $\epsilon:C^{h,k,j}\to C^{h-1,k,j}$ is bounded.
\end{pf}

\subsection{Coverings by large balls}

\begin{pro}
\label{SdeRham}
Let $X$ be a bounded geometry simplicial complex with uniformly vanishing cohomology up to degree $L$. Let $Y$ be its 0-skeleton. For every $R\in \mathbb{N}$, $R\geq 1$, consider the covering of $Y$ by balls of radius $R$, and its nerve $T^R$. The inclusion $X\subset T^R$ induces a multiplicative topological isomorphism in $\ell^\pi$-cohomology up to degree $L$, and in exact $\ell^\pi$-cohomology in degree $L+1$.
\end{pro}

\begin{pf}
Lemma \ref{leraybis} applies as in the proof of Corollary \ref{deRham}. It provides an isomorphism between cohomology at bi-degrees $(h,0)$ and $(0,h)$ for all $h$. In degrees $\leq L$, it maps cohomology of $X$ to cohomology of the nerve. In degree $L+1$, it maps exact cohomology of $X$ to exact cohomology of the nerve.
\end{pf}

\begin{rem}
Here, the cohomology isomorphism arises from a homotopy of complexes.
\end{rem}

\section{Quasi-isometry invariance}
\label{qi inv}

The above discussion suggests to use the similarity between $\ell^\pi$-cochains of a covering and Alexander-Spanier cochains, a purely metric notion. 

\subsection{Alexander-Spanier cochains}
\label{pack}

\begin{dfi}
Let $X$ be a metric space. Given $r>0$, the \emph{Rips complex of size $S$} $T_S$ of $X$ is the simplicial complex whose vertices are all points of $X$, and where a set of $k+1$ distinct vertices spans a $k$-simplex if and only if it is contained in some ball of radius $S$. Its simplicial cochains are called \emph{Alexander-Spanier cochains} of size $S$.
\end{dfi}

The following definitions are taken from \cite{lsc}.

\begin{dfi}
Let $X$ be a metric space. Given $0< R\leq S<+\infty$ and $\ell\geq 1$, a $(\ell,R,S)$-packing is a collection of balls $B_j$ such that
\begin{enumerate}
  \item the radii belong to the interval $[R,S]$,
  \item the concentric balls $\ell B_j$ are pairwise disjoint.
\end{enumerate} 
\end{dfi}

\begin{dfi}
Let $\kappa$ be an Alexander-Spanier $k$-cochain of size $S$. Its \emph{packing $\ell^p$ norm} is defined by
\begin{eqnarray*}
\|\kappa\|_{p,\ell,R,S}=\sup_{(\ell,R,S)\mathrm{-packings\,\{B_j\}}}\left(\sum_{j}\sup_{x_0,\dots, x_k\in B_j}|\kappa(x_0,\dots, x_k)|^p\right)^{1/p}.
\end{eqnarray*}
This defines a Banach space $AL^{p,k}_{\ell,R,S}(X)$.
\end{dfi}

Given $\pi=(p_0,\ldots,p_n,\ldots)$, the spaces
\begin{eqnarray*}
AL^{\pi,k}_{\ell,R,S}(X)=AL^{p_k,k}_{\ell,R,S}(X)\cap d^{-1}(AL^{p_{k+1},k+1}_{\ell,R,S}(X))
\end{eqnarray*}
form a complex of Banach spaces, whose cohomology is the \emph{packing $\ell^\pi$-cohomology} of $X$. It has a forgetful map to ordinary cohomology, whose kernel is the \emph{exact packing $\ell^\pi$-cohomology} of $X$.

\subsection{Changing size}

An Alexander-Spanier cochain of size $S\geq 1$ determines an Alexander-Spanier cochain of size 1, by restriction, whence a map $AS_{\ell,R,S}^{\cdot}(Y)\to AS_{1,1,1}^{\cdot}(Y)$, where the domain only depends on $S$ whereas parameters $\ell\geq 1$ and $R>0$ merely influence the norm. 

\begin{pro}
\label{sizech}
Let $X$ be a simplicial complex with bounded geometry and uniformly vanishing cohomology up to degree $L$. Let $Y$ be its 0-skeleton. For every integer $S$, every $0<R\leq S$ and $\ell\geq 1$, the forgetful map $AS_{\ell,R,S}^{\cdot}(Y)\to AS_{1,1,1}^{\cdot}(Y)$ induces a multiplicative topological isomorphism in $\ell^\pi$-cohomology up to degree $L$, and in exact cohomology $EH_{\ell,R,S}^{\pi,L+1}(Y)\to EH_{1,1,1}^{\pi,L+1}(Y)$ in degree $L+1$.
\end{pro}

\begin{pf}
Cochains of size 1 coincide with simplicial cochains of $X$. The counting $\ell^\pi$ norm coincides with the packing $\ell^\pi$ norm at size 1, up to a multiplicative constant depending on the local geometry of $X$. 

By construction, a collection of $S$-balls in $Y$ has a nonempty intersection if and only if their centers belong to the same $S$-ball. Thus the Rips complex of size $S$ coincides with the nerve of the covering by $S$-balls coincide with Alexander-Spanier cochains of size $S$. Let us compare norms. In nerve notation, the packing $\ell^p$-norm reads
\begin{eqnarray*}
\|\kappa\|^{p}_{p,\ell,S,S}=\sup_{J\, (\ell S)\mathrm{-separated\, subset\, of }\,Y}\sum_{j\in J}\sup_{\{i_0,\ldots,i_h\,;\,j\in U_{i_0,\ldots,i_h}\}}|\kappa(i_0,\ldots,i_h)|^p.
\end{eqnarray*}
This is always less than
\begin{eqnarray*}
\sum_{j\in Y}\sum_{\{i_0,\ldots,i_h\,;\,j\in U_{i_0,\ldots,i_h}\}}|\kappa(i_0,\ldots,i_h)|^p\leq V(S)\sum_{i_0,\ldots,i_h}\kappa(i_0,\ldots,i_h)|^p,
\end{eqnarray*}
where $V(S)$ is an upper bound for the number of vertices in an $S$-ball. Indeed, a multi-index $i_0,\ldots,i_h$ arises in the sum at most as many times as there are vertices in $U_{i_0,\ldots,i_h}$, and this is less than $V(S)$. The same crude bound remains valid for $\|\kappa\|_{p,\ell,R,S}$ for all $R\leq S$. Conversely, pick, for each $h$-simplex $i_0,\ldots,i_h$ a $j\in Y$ such that $i_0,\ldots,i_h\subset B(j,S)$, denote it by $j(i_0,\ldots,i_h)$. Assume that $Y$ can be covered with at most $N$ $(\ell S)$-separated subsets $J$. For each of them, 
\begin{eqnarray*}
\sum_{\{i_0,\ldots,i_h\,;\,j(i_0,\ldots,i_h)\in J\}}&&|\kappa(i_0,\ldots,i_h)|^p\\
&=& \sum_{j\in J}\sum_{\{i_0,\ldots,i_h\,;\,j(i_0,\ldots,i_h)=j\}}|\kappa(i_0,\ldots,i_h)|^p\\
&\leq&
V(S)\sum_{j\in J}\sup_{\{i_0,\ldots,i_h\,;\,j\in U_{i_0,\ldots,i_h}\}}|\kappa(i_0,\ldots,i_h)|^p,
\end{eqnarray*}
hence
\begin{eqnarray*}
\sum_{i_0,\ldots,i_h}|\kappa(i_0,\ldots,i_h)|^p \leq
N V(S)\|\kappa\|^{p}_{p,\ell,S,S}\leq
N V(S)\|\kappa\|^{p}_{p,\ell,R,S}.
\end{eqnarray*}
To get an upper bound on $N$, let us construct inductively a colouring of $Y$ with values in $\{0,\ldots,V(\ell S)\}$. Pick an origin $o$, and colour it 0. Assume a finite part $A$ of $Y$ has already been coloured, pick a point $y$ among the uncoloured points which are closest to $o$, choose its colour among those which are not already used in $B(y,\ell S)\cap A$. This is possible since $|B(y,\ell S)|<V(\ell S)+1$. In such a way, one colours all of $Y$, and each set $J$ of points of equal colour is $(\ell S)$-separated. So $N=V(\ell S)+1$ is appropriate.

This shows that the counting $\ell^\pi$ norm on cochains of the covering and the packing norm are equivalent, with constants depending only on the geometry of $X$ at scale $S$, i.e. on $S$ only. Thus $\ell^\pi$-cohomology of the nerve $T^{S}$ coincides with packing $\ell^\pi$-cohomology at size $S$, with equivalent norms. The inclusion of nerves corresponds to the forgetful map for cochains. Thus the statement is a reformulation of Proposition \ref{SdeRham}.
\end{pf}

\subsection{Invariance}

Say a map $f:X\to X'$ between metric spaces is a \emph{coarse embedding} if for every $T>0$, there exists $T'(T)>0$ such that for every $T$-ball $B$ of $X$ and every $T$-ball $B'$ of $X'$, $f(B)$ and $f^{-1}(B')$ are contained in $T'$-balls. A \emph{quasi-isometry} is a pair of coarse embeddings $f:X\to X'$ and $g:X'\to X$ such that $f\circ g$ and $g\circ f$ are a bounded distance away from identity.

Packing cohomology is natural under coarse embeddings, up to a loss on size. Furthermore, embeddings which are a bounded distance away from each other induce the same morphism in packing cohomology.

\begin{pro}
\label{qii}
Let $f:X\to X'$ be a coarse embedding between metric spaces. Then for every $R>0$, $R\leq S<\infty$ and $\ell'\geq 1$, there exist $R'>0$, $S'<\infty$ and $\ell\geq 1$ such that $f$ induces a multiplicative morphism $f^*:L^{\pi}_{\ell',R',S'}H^{\cdot}(X')\to L^{\pi}_{\ell,R,S}H^{\cdot}(X)$.

If $g:X\to X'$ satisfies $\sup_{x\in X}d(f(x),g(x))<+\infty$, then $g$ is a coarse embedding as well, and $f^*=g^*$.
\end{pro}

\begin{pf}
Given a size $r>0$, by definition of a coarse embedding, there exists $r'(r)$ such that composition with $f$ maps cochains of size $r'$ to cochains of size at least $r$.

Given $0<R\leq S<+\infty$ and $\ell'\geq 1$, let $R'=S'=r'(S)$, let $T=2\ell'S'$ and $\ell=r'(T)/R$. Then $f$ maps $(\ell,R,S)$-packings to $(\ell',R',S')$-packings. Thus composition with $f$ is bounded in suitable packing norms. It commutes with $d$ and with cup-product. Therefore it induces a multiplicative morphism $f^*:L^{\pi}_{\ell',R',S'}H^{\cdot}(X')\to L^{\pi}_{\ell,R,S}H^{\cdot}(X)$.

Given simplices $\Delta=\{x_{0},\ldots,x_{k}\}$ and $\Delta'=\{x'_{0},\ldots,x'_{k}\}$ of $T_S(X)$, the {\sl prism} $b(\Delta,\Delta')$, obtained by triangulating the product of a simplex and an interval, is defined by
$$
b(x_{0},\ldots,x_{k};x'_{0},\ldots,x'_{k})=
\sum_{i=0}^{k} (-1)^{i}(x_{0},\ldots,x_{i-1},x_{i},x'_{i},x'_{i+1},\ldots,x'_{k}).
$$
It satisfies
$$
\partial b(\Delta,\Delta')=\Delta'-\Delta-\sum_{j=0}^{k} 
(-1)^{j}b(\partial_{j}\Delta,\partial_{j}\Delta').
$$
Assume that $\sup_{x\in X}d(f(x),g(x))\leq\epsilon$. If $\kappa$ is a $k$-cochain of size $S+\epsilon$ and $\Delta$ a simplex of size $S$, set $(B\kappa)(\Delta)=\kappa(b(f(\Delta),g(\Delta)))$. Then
\begin{eqnarray*}
dB+Bd=\kappa\circ g-\kappa\circ f.
\end{eqnarray*}
For all $\ell\geq 1$, 
\begin{eqnarray*}
\|B\|_{L^{p}_{\ell,R-\epsilon,S+\epsilon}\to L^{p}_{\ell,R,S}}\leq (k+1)^{1/p}.
\end{eqnarray*}
This shows that $f^*=g^*$ on cochains of sufficiently large size.
\end{pf}

\subsection{Packing $\ell^\pi$-cohomology equals $L^\pi$-cohomology}

Let $X$ be a bounded geometry Riemannian manifold. Pick a uniform sequence of nested coverings $\mathcal{U}^0,\ldots,\mathcal{U}^L$. Up to rescaling once and for all the metric on $X$, one can assume that coverings have the following properties:
\begin{enumerate}
  \item Each $U_i^0$ contains a unit ball $B_i$, and these balls are disjoint.
  \item Each $U_i^0$ is contractible in $U_i^1$.
  \item The diameters of $U_i^1$ are bounded.
\end{enumerate}
Under these assumptions, the 0-skeleton $Y$ of the nerve $T$ of the covering $\mathcal{U}^0$ is quasi-isometric to $X$. Indeed, the map that sends vertices to centers of balls $B_i$ is bi-Lipschitz, its image is $D$-dense for some finite $D$. The Rips complex of size 1 of $Y$ coincides with $T_0$, so packing $\ell^\pi$-cohomology of $Y$ at size 1 coincides with simplicial $\ell^\pi$-cohomology of $T^0$ (packing norms for a uniformly discrete metric spaces are equivalent to counting norms) at size 1. According to Corollary \ref{deRham}, this is equal to de Rham $L^\pi$-cohomology of $X$.

We can now proceed to the proof of Theorem \ref{qib}. Since the covering pieces are contractible in unions of boundedly many pieces), the natural map of $X$ to the nerve, given by a partition of unity, is a homotopy equivalence. Therefore uniform vanishing of cohomology passes from $X$ to the nerve.

By Corollary \ref{deRham}, $L^\pi$-cohomology of $X$ is isomorphic to $\ell^\pi$-cohomology of the nerve, which in turn coincides with packing $\ell^\pi$-cohomology of $Y$ at all sizes by Proposition \ref{sizech}.

The inclusion $i:Y\to X$ is a quasi-isometry. According to Proposition \ref{qii}, for all $0<R\leq S<\infty$ and all $\ell''$, there exist $\ell,\ell'\geq 1$ and $R',R''>0$, $S',S''<\infty$ such that $i$ and its inverse induce maps up to degree $L$
\begin{eqnarray*}
L^{\pi}_{\ell'',R'',S''}H^{\cdot}(Y)\to L^{\pi}_{\ell',R',S'}H^{\cdot}(X)\to L^{\pi}_{\ell,R,S}H^{\cdot}(Y),
\end{eqnarray*}
and in the reverse direction. The composition coincides with the forgetful map (Proposition \ref{qii}), which is an isomorphism, by Proposition \ref{sizech}. Therefore 
$$
i^*:L^{\pi}_{\ell',R',S'}H^{\cdot}(X)\to L^{\pi}_{\ell,R,S}H^{\cdot}(Y)
$$
is an isomorphism. This proves that de Rham and packing $L^\pi$-cohomo\-logies of $X$ are isomorphic up to degree $L$, and that forgetful maps induce isomorphisms in the packing $\ell^\pi$-cohomology of $X$ up to degree $L$. In degree $L+1$, the result persists provided one considers exact $L^\pi$-cohomology.

A quasi-isometry between manifolds $f:X\to X'$ gives rise to cohomology maps in both directions with a loss on size, whose compositions coincide with forgetful maps. Since forgetful maps are isomorphisms, $f^*$ is an isomorphism up to degree $L$, and an isomorphism on exact $L^\pi$-cohomology in degree $L+1$. The isomorphism is topological unless $(p,q)=(1,\frac{n}{n-1})$, or $(n,\infty)$, as observed in Remark \ref{specialcases}.

\subsection{$(L^q,L^p +L^r)$-cohomology}

...TO DO...

\section{Contact manifolds}
\label{contact}

\subsection{Sub-Riemannian contact manifolds}

A \emph{sub-Riemannian manifold} is the data of a manifold $M$, a smooth sub-bundle $H\subset TM$, and a smooth field of Euclidean structures on $H$.

A smooth codimension 1 sub-bundle $H\subset TM$ can be defined as the kernel of a smooth 1-form $\theta$. Up to a scale, the restriction of $d\theta$ to $H$ does not depend on the choice of $\theta$. Say $(M,H)$ is a \emph{contact manifold} if $d\theta_{|H}$ is non-degenerate.

A sub-Riemannian metric on a $2m+1$-dimensional contact manifold extends canonically into a Riemannian metric. Indeed, there is a unique contact form $\theta$ such that $\frac{1}{m!}{(d\theta)^m}_{|H}$ equals the Euclidean volume form on $H$. This contact form is smooth, the kernel of $d\theta$ defines a complement to $H$ carrying the Reeb vectorfield $\rho$, normalized so that $\langle\theta,\rho\rangle=1$, hence the unique Riemannian metric which makes $\rho\perp H$ and $|\rho|^2=1$.

\begin{rem}
A sub-Riemannian contact manifold has bounded geometry (see Definition \ref{defbddgeocontact}) if and only if the corresponding Riemannian metric has bounded geometry. 
\end{rem}

\subsection{Rumin's complex}

On a $2m+1$-dimensional contact manifold, consider the algebra $\Omega^\cdot$ of smooth differential forms, let $\mathcal{I}^\cdot$ denote the ideal generated by $1$-forms vanishing on $H$, let $\mathcal{J}^\cdot$ denote its annihilator. The exterior differential descends (resp. restricts) to an operator $d_c:\Omega^\cdot/\mathcal{I}^\cdot \to \Omega^\cdot/\mathcal{I}^\cdot$ (resp. $d_c:\mathcal{J}^\cdot\to \mathcal{J}^\cdot$). Note that $\mathcal{I}^h=0$ for $h\geq m+1$ and $\mathcal{J}^h=0$ for $h\leq m$. In \cite{Rumin1}, M. Rumin defines a second order linear differential operator $d_c:\Omega^m/\mathcal{I}^m \to \mathcal{J}^{m+1}$ which connects $\Omega^\cdot/\mathcal{I}^\cdot$ and $\mathcal{J}^\cdot$ into a complex ($d_c\circ d_c=0$) which can be used to compute cohomology. $\Omega^\cdot/\mathcal{I}^\cdot$ and $\mathcal{J}^\cdot$ identify with spaces of smooth sections of bundles $E_0^h$, $h=0,\ldots,2m+1$, which inherit Euclidean structures, therefore $L^p$ norms make sense.

\begin{thm}
\label{qir}
Assume that $1\leq p\leq q\leq \infty$ and $\frac{1}{p}-\frac{1}{q}\leq \frac{1}{2m+2}$ (to be replaced with $\frac{1}{p}-\frac{1}{q}\leq \frac{2}{2m+2}$ when degree $m+1$ cohomology is considered). 

Consider the class of contact sub-Riemannian manifolds with the following properties.
\begin{enumerate}
\item Dimension equals $2m+1$.
  \item Bounded geometry.
  \item Uniform vanishing of cohomology up to degree $k-1$.
\end{enumerate}
If $p=2m+2$, $q=\infty$ and $k=1$, one should replace $L^{\infty,p}$-cohomology with $L^{BMO,p}$-cohomology.

Assume that $k\leq m$. For $X$ in this class, and up to degree $k-1$, the $L^{q,p}$-cohomology of Rumin's complex and the packing $\ell^{q,p}$-cohomology of $X$ at all sizes are isomorphic as vectorspaces. In degree $k$, it is the exact $L^{q,p}$-cohomology of Rumin's complex which is isomorphic to packing $\ell^{q,p}$-cohomology.

If $k\geq m+1$, the same conclusion holds in non-limiting cases, i.e. if either $p>1$, $q<\infty$ or $\frac{1}{p}-\frac{1}{q}< \frac{1}{2m+2}$ (resp. $\frac{1}{p}-\frac{1}{q}< \frac{2}{2m+2}$ in degree $m+1$). 
\end{thm}

The given sub-Riemannian metric and the corresponding Riemannian metric are quasi-isometric, so their packing $\ell^\pi$-cohomologies are isomorphic. Therefore, under the assumptions of Theorem \ref{qir}, the Rumin complex can be used to compute packing $\ell^\pi$-cohomology.

\begin{exa}\label{cohohm}
If $k\in \{0,\ldots,2m+1\}$, $1<p\leq q<\infty$ satisfy
\begin{eqnarray*}
\frac{1}{p}-\frac{1}{q}=\begin{cases}
\frac{1}{2m+2}      & \text{ if }k\not=m+1, \\
\frac{2}{2m+2}      & \text{ if }k=m+1,
\end{cases}
\end{eqnarray*}
then $H^{q,p,k}(\he{m})=0$.
\end{exa}
Indeed, the existence of global homotopy operators, Theorem \ref{homotopy formulas} of \cite{SmP}, implies that the cohomology of the Rumin complex vanishes, and thus $\ell^{q,p}$ cohomology vanishes by Theorem \ref{qir}.

\subsection{Cutting-off Rumin differential forms}

The proof of Theorem \ref{qir} follows the same lines as Theorem \ref{qib}. The local model for $2m+1$-dimensional sub-Riemannian contact manifolds is the Heisenberg group $\H^m$ equipped with its left-invariant contact structure and a left-invariant Euclidean structure on it. The local ingredients are
\begin{enumerate}
  \item An inverse of the analytic differential $d$ on balls, possibly with a loss on the domain: this is given by Poincar\'e inequalities. According to \cite{SmP}, Poincar\'e inequalities are valid in balls of $\H^m$ with respect to Rumin's differentials $d_c$. The fact that Rumin's differential in degree $m$ is second order allows the broader inequality $\frac{1}{p}-\frac{1}{q}\leq\frac{2}{2m+2}$ in degree $m+1$.
  \item An inverse of the combinatorial coboundary $\delta$. 
\end{enumerate}
In Lemma \ref{pq}, the following inverse $\epsilon$ was used,
\begin{eqnarray*}
\epsilon(\phi)_{i_0 \ldots i_{h-1}}=\sum_{j}\chi_j \phi_{j i_0 \ldots i_{h-1}}.
\end{eqnarray*}
It is bounded on $L^p$. One needs it to be bounded on $L^{q}\cap d_c^{-1}L^p$. In Lemma \ref{pq}, this relies on Leibniz' formula
\begin{eqnarray*}
d(\zeta\alpha)=d\zeta\wedge\alpha+\zeta d\alpha.
\end{eqnarray*}
A difficulty arises in the contact case since the middle $d_c$ is second order: Leibniz formula reads
\begin{eqnarray*}
d_c(\zeta\alpha)=\zeta d_c(\alpha)+P(\nabla\zeta,\nabla\alpha)+Q(\nabla^2\zeta,\alpha),
\end{eqnarray*}
where $Q$ does not differentiate $\alpha$, so it is bounded on $L^q$, but $P$ does depend on all horizontal first derivatives of $\alpha$, and is not expressible in terms of $d_c\alpha$ only. The solution consists in passing to a homotopy equivalent complex of forms whose horizontal first derivatives are controlled. This modification, needed only to handle degrees $\geq m+1$, does not work in limiting cases yet.

\subsection{The $W^{1,\pi}$-Rumin complex}

\begin{dfi}
Let $M$ be a sub-Riemannian contact manifold. Let $W^{1,p}(M,E_0^k)$ denote the space of degree $k$ Rumin forms which satisfy, in the sense of distributions,
\begin{eqnarray*}
|\alpha|\in L^p, \quad|\nabla \alpha|\in L^p.
\end{eqnarray*}
Given a vector of exponents $\pi=(p_0,\ldots,p_{2m+1})$, let 
\begin{eqnarray*}
\Omega_c^{\pi,k}(M)&=&L^{p_k}(M,E_0^k)\cap d_c^{-1}(L^{p_{k+1}}(M,E_0^{k+1})),\\
\Omega_W^{\pi,k}(M)&=&W^{1,p_k}(M,E_0^k)\cap d_c^{-1}(W^{1,p_{k+1}}(M,E_0^{k+1}))
\end{eqnarray*}
denote the two complexes one can form with Rumin forms: the $L^\pi$
Rumin complex and the $W^{1,\pi}$ Rumin complex.
\end{dfi}

Multiplication with a smooth function maps $\Omega_W^{\pi,\cdot}$ to $\Omega_c^{\pi,\cdot}$.

Here is the relevant Poincar\'e inequality. It is valid provided the following inequalities hold.

\begin{eqnarray}\label{kpq}
1 \leq p\leq q \leq\infty,\quad \frac{1}{p}-\frac{1}{q}\leq\begin{cases}
\frac{1}{2m+2}      & \text{ if }k\not=m+1, \\
\frac{2}{2m+2}      & \text{ if }k=m+1,
\end{cases}
\end{eqnarray}
We speak of a \emph{limiting case} when $\frac{1}{p}-\frac{1}{q}=\frac{1}{2m+2}$ (resp. $=\frac{2}{2m+2}$ in degree $m+1$) and one of $p$ and $q$ equals 1 (resp. $\infty$). 

\begin{lem}[\ref{improved Poincare} of \cite{SmP}]\label{poincaredc}
Assume that $(k,p,q)$ satisfy inequations (\ref{kpq}) above. There exists $\lambda>1$ and $C(\lambda)$ such that the following holds. Let $B=B(e,1)$ and $B''=B(e,\lambda)$ be concentric balls of $\he{n}$. 

Assume first that $(p,q,k)\notin\{(1,\frac{2m+2}{2m+1},2m+1),(2m+2,\infty,1)\}$. For every closed differential $k$-form $\omega$ on $B''$, there exists a differential $k-1$-form $\phi$ on $B$ such that $d\phi=\omega_{|B}$ and
\begin{eqnarray*}
\|\phi\|_{L^{q}(B)}\leq C\,\|\omega\|_{L^{p}(B'')}.\quad \quad \quad (\mathbb{H}-Poincare_{q,p}(k))
\end{eqnarray*}

If $p=1$, $q=\frac{2m+2}{2m+1}$ and $k=2m+1$, inequality $(\mathbb{H}-Poincare_{q,p}(k))$ is replaced with
\begin{eqnarray*}
\|\phi\|_{L^{\frac{2m+2}{2m+1}}(B)}\leq C\,\|\omega\|_{\mathcal{H}^{1}(B'')}.\quad \quad \quad (\mathbb{H}-Poincare_{\frac{2m+2}{2m+1},\mathcal{H}^1}(2m+1))
\end{eqnarray*}

If $p=2m+2$, $q=\infty$ and $k=1$, inequality $(Poincare_{q,p}(k))$ is replaced with
\begin{eqnarray*}
\|\phi\|_{BMO(B)}\leq C\,\|\omega\|_{L^{2m+2}(B'')}.\quad \quad \quad (\mathbb{H}-Poincare_{BMO,2m+2}(1))
\end{eqnarray*}

In non-limiting cases, for every closed $h$-form $\alpha\in W^{1,p}(B'',E_0^h)$, there exists an $h-1$-form $\beta\in W^{1,q}(B,E_0^{h-1})$ such that $d_c\beta=\alpha_{|B}$, and
\begin{eqnarray*}
\|\beta\|_{W^{1,q}(B)}\leq C\,\|\alpha\|_{W^{1,p}(B'')}.
\end{eqnarray*}
\end{lem}

\subsection{Back to the $L^p$ Rumin complex}

To overcome the fact that the inverse $\epsilon$ of the \v Cech coboundary $\delta$ involves a loss of differentiability (it merely maps $\Omega_W^{1,\pi}$ to $\Omega_c^{1,\pi}$), we shall use a local smoothing procedure, provided again by \cite{SmP}. 

\begin{lem}[\ref{smoothing} of \cite{SmP}]\label{smoothingdc}
Let $B=B(e,1)$ and $B'=B(e,2)$ be concentric balls of $\he{m}$. There exist operators $S$ and $T$ from smooth forms on $B'$ to smooth forms on $B$ which satisfy $S+d_c T+Td_c=R_B$, the restriction of forms to $B$. For every $(k,p,q)$ satisfying inequations (\ref{kpq}) above, excluding limiting cases, and for every $L\in\mathbb{N}$, these operators extend to bounded operators $T:L^p(B',E_0^\cdot)\to L^{p}(B,E_0^{\cdot -1})$ and $S:L^p(B',E_0^\cdot)\to W^{L,p}(B,E_0^{\cdot})$. Furthermore, $T$ is bounded $W^{1,p}(B',E_0^\cdot)\to W^{1,q}(B,E_0^{\cdot -1})$. 
\end{lem}

Since the smoothing operator $S$ is only locally defined, it does not directly provide us with a homotopy equivalence $\Omega_{c}^{\cdot,\pi}\to \Omega_{W}^{\cdot,\pi}$. We must pass via the bi-complexes $C_c^{\cdot,\cdot,j}$, $j=0,\ldots,L$, associated to a uniform sequence of nested coverings.

Proposition \ref{existcover} allows to adjust the ratio of radii of the model Heisenberg balls $Z'\subset Z$. Choose this ratio to be $\geq 2$, in order that Lemma \ref{smoothingdc} be applicable and yields operators $r\circ S$ and $r\circ T$ defined on the Rumin bi-complex
\begin{eqnarray*}
rS,\,rT:C_c^{\cdot,\cdot,j}\to C_c^{\cdot,\cdot,j-1}
\end{eqnarray*}
constructed from the $L^\pi$ Rumin bicomplex. Here, $d'=\delta$ is the \v Cech coboundary, and $d''=(-1)^h d_c$ is the Rumin differential (up to sign). Let $\epsilon$ denote the operator defined in Lemma \ref{pq}, which satisfies
$\epsilon\delta+\delta\epsilon=1$. Let us compute
\begin{eqnarray*}
(d'+d'')(rS\epsilon+(-1)^h rT)
&=&d'rS\epsilon-rSd'\epsilon+rSd'\epsilon+d''rS\epsilon\\
&&+(-1)^h d''rT+(-1)^h d''rT,\\
(rS\epsilon+(-1)^h rT)(d'+d'')
&=&rS\epsilon d'+rS\epsilon d''\\
&&+(-1)^h rTd''+(-1)^h rTd'.
\end{eqnarray*}
Since $d'\epsilon+\epsilon d'=1$ on $C^{h,k,j}$ and $rS+(-1)^h d''rT+(-1)^h rTd''=r$,
\begin{eqnarray*}
(d'+d'')(rS\epsilon+(-1)^h rT)+(rS\epsilon+(-1)^h rT)(d'+d'')
&=&r+U+V,
\end{eqnarray*}
where 
\begin{eqnarray*}
U=d'rS\epsilon-rSd'\epsilon+rS\epsilon d''+d''rS\epsilon
\end{eqnarray*}
is smoothing, and 
\begin{eqnarray*}
V=(-1)^h (d'rT+rTd')
\end{eqnarray*}
has bi-degree $(1,-1)$. Denote by $D=d'+d''$, $B=rS\epsilon+(-1)^h rT$ and $W=U+V$. Note that $WD=DW$. One can iterate identity $DB+BD=r+W$ as follows. Write
\begin{eqnarray*}
DBW+BWD&=&(DB+BD)W=rW+W^2,\\
DrB+rBW&=& r(DB+BD)=r^2+rW
\end{eqnarray*}
and substract,
\begin{eqnarray*}
D(rB-BW)+(rB-BW)B=r^2-W^2.
\end{eqnarray*}
Ultimately, we find a polynomial $P$ in $r$ and $W$ such that $DBP+BPD=r^L-(-1)^L W^L$. Since $V$ has bi-degree $(1,-1)$, $V^L=0$, hence $W^L$ is a sum of words in $U$ and $V$ such that each term has at least a $U$ in it, hence is smoothing. This provides a homotopy of $r^L$ to a bounded operator $C_c^{\cdot,\cdot,L}\to C_W^{\cdot,\cdot,0}$. 

Up to the cost of enlarging the number of nested coverings required to $L^2$, we can follow each use of $\epsilon$ with a use of $W^L$, and return to the bi-complexes $C_W^{\cdot,\cdot,j}$ without changing homotopy types. This makes it possible to apply Lemma \ref{leraybis} as in the proof of Corollary \ref{deRham}. This proves Theorem \ref{qir}.

\bigskip

\noindent
Pierre Pansu 
\par\noindent Laboratoire de Math\'ematiques d'Orsay,
\par\noindent Univ. Paris-Sud, CNRS, Universit\'e
Paris-Saclay
\par\noindent 91405 Orsay, France.
\par\noindent
e-mail: pierre.pansu@math.u-psud.fr

\end{document}